\documentclass[a4paper,11pt,twoside]{amsart}

\usepackage{amsmath}
\usepackage{amsthm}
\usepackage{amssymb}
\usepackage{hyperref}

\newtheorem{thm}{Theorem}[section]
\newtheorem{prop}[thm]{Proposition}
\newtheorem{lem}[thm]{Lemma}
\newtheorem{cor}[thm]{Corollary}

\theoremstyle{definition}
\newtheorem{definition}[thm]{Definition}
\newtheorem{remark}[thm]{Remark}
\newtheorem{ex}[thm]{Example}

\newcommand{\supp}{\operatorname{supp}}
\newcommand{\Sing}{\operatorname{Sing}}

\begin{document}

\title{Singular cubic fourfolds containing a plane}
\author{Paolo Stellari}

\address{Dipartimento di Matematica ``F.
Enriques'', Universit{\`a} degli Studi di Milano, Via Cesare
Saldini 50, 20133 Milano, Italy}

\email{paolo.stellari@unimi.it}

\keywords{Cubic 4-folds, plane sextics and odd
theta-characteristics}

\subjclass[2000]{14J35, 14H50}

\begin{abstract}
In this paper we consider cubic 4-folds containing a plane whose
discriminant curve is a reduced nodal plane sextic. In particular,
we describe the singular points of such cubic 4-folds and we give
an estimate of the rank of the free abelian group generated by the
equivalence classes of the algebraic cycles of codimension 2.
Moreover, we describe some conditions on the geometry of the plane
sextics so that all the associated cubic 4-folds are singular and
we construct a family of smooth rational cubic 4-folds whose
discriminant curve is reduced but reducible.
\end{abstract}

\maketitle

\section{Introduction}


In this paper, a cubic fourfold is a hypersurface of degree 3 in
$\mathbb{P}^5$. Beauville and Voisin (see \cite{B1}, \cite{B2} and
\cite{V}) proved that if $C$ is a smooth plane sextic and $\theta$
is an odd theta-characteristic such that $h^0(C,\theta)=1$, then
the pair $(C,\theta)$ determines a cubic 4-fold containing a
plane. In \cite{CMF}, Friedman and Casalaina-Martin showed that
when $C$ is an irreducible nodal plane quintic and $\theta$ is the
push-forward of a theta-characteristic on the complete
normalization $\tilde C$ of $C$, then the corresponding cubic
3-fold is smooth.

The aim of this paper is to study the relation between cubic
fourfolds containing a plane and reduced nodal plane sextics
endowed with odd theta-characteristics. Our main result is Theorem
\ref{thm:3.1} which gives an explicit description of the singular
points of a cubic 4-fold $X$ containing a plane $P$ which is
associated to a reduced nodal plane sextic $C$ and to a
generalized theta-characteristic $\theta$ on $C$ with
$h^0(C,\theta)=1$. In particular, this result analyzes the
connections with the geometry of $\theta$ and $C$, it counts the
number of singular points and it gives an estimate of the rank of
the free abelian group $\mathrm{NS}_2(X)$ generated by the
equivalence classes of the algebraic cycles of codimension 2 in
$X$. It is worth noticing that, due to the work of Hassett
(\cite{H1} and \cite{H2}), a better understanding of the
lattice-theoretic properties of this group can give interesting
information about the geometry of the cubic 4-folds (for example,
their rationality).

Proposition \ref{thm:4.1} gives some sufficient conditions on the
geometry of a plane sextic which imply that all the associated
cubic 4-folds are singular. Proposition \ref{thm:4.4} shows the
existence of a family of smooth rational cubic 4-folds whose
discriminant curve is a reduced nodal plane sextic.

In Section \ref{sec:thetacubic} we prove some general facts about
cubic 4-folds and theta-characteristics for nodal reduced (but
possibly reducible) curves. We will work over the complex numbers.


\section{Cubic 4-folds and
theta-characteristics}\label{sec:thetacubic}


Given a smooth curve $C$, a {\it theta-characteristic} is a line
bundle $\theta$ such that
\[
\theta^{\otimes 2}=\omega_C,
\]
where $\omega_C$ is the canonical sheaf of $C$. We say that a
theta-characteristic is {\it even} if $h^0(C,\theta)\equiv 0
\pmod{2}$ and {\it odd} if $h^0(C,\theta)\equiv 1 \pmod{2}$. In
particular, by the adjunction formula, if $C$ is a smooth plane
sextic, a theta-characteristic on $C$ is a line bundle such that
\[
\theta^{\otimes 2}=\mathcal{O}_{\mathbb{P}^2}(3)|_C.
\]
It is a classical result (see for example \cite{C}) that the
number of non-isomorphic theta-characteristics is equal to
$2^{2g}$, where $g$ is the genus of $C$. More precisely, there are
\mbox{$2^{g-1}(2^g+1)$} even and \mbox{$2^{g-1}(2^g-1)$} odd
theta-characteristics.

When $C$ is a singular curve, the general picture becomes slightly
different. We recall the following definition:


\begin{definition}\label{def:2.1}
A \emph{stable spin curve} is a pair $(Y,L)$, where $Y$ is a reduced connected curve with the
following properties:

(a.1) it has only ordinary double points;

(b.1) each smooth rational component contains at least 2 nodes;

(c.1) two smooth rational components never meet each other

\noindent (we will often call {\it exceptional component} a smooth
rational component of $C$ satisfying property (b.1)). Moreover $L$
is a line bundle such that, if $Z=\overline{Y-(\cup_{i\in I}
E_i)}$, where $\{E_i:i\in I\}$ is the set of the exceptional
components of $Y$, then

(a.2) $(L|_Z)^{\otimes 2}\cong \omega_Z$;

(b.2) $L|_{E_i}\cong\mathcal{O}_{E_i}(1)$, for all $i\in I$.

A \emph{stable model} of a stable spin curve $(Y,L)$ is a curve
$C$ which is obtained by contracting all the exceptional
components of $Y$. The map $\nu:Y\rightarrow C$ contracting the
exceptional components of $Y$ is the {\it contraction
map}.\end{definition}


Cornalba showed (see \cite{C}, page 566) that the isomorphism
classes of stable spin curves whose stable model is a nodal curve
$C$ are the natural analogues of the theta-characteristics for
$C$. In particular, one can prove (see \cite{CC}) that any nodal
curve is the stable model of $2^{2g}$ isomorphism classes of spin
curves (counted with multiplicity), where $g$ is the arithmetic
genus of the curve $C$.


\begin{definition}\label{def:2.2} A {\it theta-characteristic} on
a reduced nodal curve $C$ is a sheaf $\theta$ on $C$, such that
there exists a stable spin curve $(Y,L)$ whose stable model is $C$
and satisfying the following condition:
\[
\theta=\nu_* L,
\]
where $\nu: Y\rightarrow C$ is the contraction map.
\end{definition}


Also in this case, we say that a theta-characteristic $\theta$ is
even (resp. odd) if \mbox{$h^0(C,\theta)\equiv 0\pmod{2}$} (resp.
$h^0(C,\theta)\equiv 1\pmod{2}$). In particular, one can prove
that there are $2^{g-1}(2^g+1)$ even and $2^{g-1}(2^g-1)$ odd
theta-characteristics, where these sheaves are counted with
multiplicity and $g$ is the arithmetic genus of $C$.

Let $C$ be a reduced nodal curve and $\theta$ a
theta-characteristic on $C$. We can associate to $\theta$ a set of
points $S_\theta\subset\Sing (C)$ defined in the following way:
\[
S_\theta:=\{p\in\Sing(C):\theta_p\not\cong\mathcal{O}_{C,p}\}.
\]

The following lemma gives a different characterization of the set
$S_\theta$.


\begin{lem}\label{lem:2.1} Let $C$ be a reduced nodal curve
and $\theta$ a theta-characteristic such that
\mbox{$\theta=\nu_*L$}, where $(Y,L)$ is a stable spin curve whose
stable model is $C$ and $\nu$ is the contraction map. Then
\[
S_\theta=\{p\in\Sing(C):\nu^{-1}(p)\mbox{ is an exceptional
component}\}.
\]
\end{lem}


\begin{proof} By Definition \ref{def:2.2}, if $p$ is such that
$\nu^{-1}(p)$ is an  exceptional component then $\theta$ is not
invertible in $p$ and thus $p\in S_\theta$. Let $p\in S_\theta$.
Then we have only two possibilities: (1) $\nu^{-1}(p)$ is one
point or (2) $\nu^{-1}(p)$ is an exceptional component. It is easy
to see that (1) cannot occur. Indeed, in this case, by Definitions
\ref{def:2.1} and \ref{def:2.2}, $\theta$ would be locally
invertible in $p$.

If $p$ is as in case (2) then
$\theta_p/m_p\theta_p\cong{\mathbb{C}}^2$, where $m_p$ is the
maximal ideal of $\theta_p$. Indeed, let $U$ be an open
neighborhood of $p$. Then we have the following restriction map:
\[
\rho:\theta(U)\rightarrow L(U_1)\oplus L(U_2)\oplus L(E),
\]
where $U_1$ and $U_2$ are open neighborhoods of the two intersection
points $p_1$ and $p_2$ of the strict transform of $C$ with the
exceptional component $E$. Since
$\theta(U)=(\nu_*L)(U)=L(\nu^{-1}(U))$ and $L(U_1)\oplus
L(U_2)\oplus L(E)=L(U_1)\oplus L(U_2)\oplus
H^0({\mathbb{P}}^1,\mathcal{O}_{{\mathbb{P}}^1}(1))$, let
\[
\rho:s\longmapsto (s_1,s_2,v)\in L(U_1)\oplus L(U_2)\oplus
H^0({\mathbb{P}}^1,\mathcal{O}_{{\mathbb{P}}^1}(1)).
\]
We want to show that $v$ is uniquely determined by $s_1$ and $s_2$
and that every $v$ occurs in the image of $\rho$ (i.e. the values
of $s_1(p)$ and $s_2(p)$ are arbitrary). This would imply
$\theta_p/m_p\theta_p\cong{\mathbb{C}}^2$.

The fact that $v$ is determined by $s_1$ and $s_2$ follows from
easy calculations (recalling that $v$ must be properly glued with
$s_1$ and $s_2$ in $p_1$ and $p_2$). Consider the exact sequence:
\[
0\longrightarrow\mathcal{O}_{{\mathbb{P}}^1}(1)(-p_1-p_2)\longrightarrow\mathcal{O}_{{\mathbb{P}}^1}(1)\longrightarrow\mathcal{O}_{p_1}\oplus\mathcal{O}_{p_2}\longrightarrow
0,
\]
where $\mathcal{O}_{p_i}$ is the skyscraper sheaf in $p_i$. The
first sheaf is isomorphic to the sheaf
$\mathcal{O}_{{\mathbb{P}}^1}(-1)$. Thus, considering the long
exact sequence in cohomology, we get
\[
0\longrightarrow 0\longrightarrow
H^0({\mathbb{P}}^1,\mathcal{O}_{{\mathbb{P}}^1}(1))\longrightarrow
{\mathbb{C}}\oplus{\mathbb{C}}\longrightarrow 0.
\]
The isomorphism
$H^0({\mathbb{P}}^1,\mathcal{O}_{{\mathbb{P}}^1}(1))\rightarrow
{\mathbb{C}}\oplus{\mathbb{C}}$ proves that $v$ is uniquely
determined by $s_1$ and $s_2$.\end{proof}


\begin{cor}\label{cor:2.2} Let $C$ be a reduced nodal
curve and $\theta$ a theta-characteristic such that
$\theta=\nu_*L$, where $(Y,L)$ is a stable spin curve whose stable
model is $C$ and $\nu$ is the contraction map. Let
$Z:=\overline{Y-(\cup_{i\in I} E_i)}$, where $\{E_i:i\in I\}$ is
the set of the irreducible components of $C$ and $\nu_Z:=\nu|_Z$.
Then $\theta=(\nu_Z)_*(L|_Z)$.
\end{cor}


\begin{proof} This easily follows form the previous lemma and
remarks. \end{proof}


\begin{cor}\label{cor:2.3} Let $C$ be a reduced nodal
curve and $\theta$ a theta-characteristic. Let us suppose that
there are two stable spin curves $(Y_1,L_1)$ and $(Y_2,L_2)$ such
that $\theta=(\nu_1)_*L_1=(\nu_2)_*L_2$, where $\nu_1$ and $\nu_2$
are the corresponding contraction maps. Then $Y_1=Y_2$.
\end{cor}


\begin{proof} Let $E:=\{E_1,\ldots,E_m\}$ and $E':=\{E'_1,\ldots,E'_n\}$,
be the sets of the irreducible components of $Y_1$ and $Y_2$
respectively. Let $S_1$ and $S_2$ be the subsets of $\Sing(C)$
defined in the following way:
\[
S_1:=\nu_1(E)\;\;\;\;\;\;\;\;\;\;\;S_2:=\nu_2(E').
\]

Let us suppose that $S_1\not\subseteq S_2$. If $p\in S_1-S_2$, then,
$\theta=(\nu_1)_*L_1$ is not invertible in $p$. We have also
$\theta=(\nu_2)_*L_2$ and, due to Lemma \ref{lem:2.1}, $\theta$ must
be invertible in $p$. This gives a contradiction. The same holds if
$S_2\not\subseteq S_1$. Hence $S_1=S_2$.\end{proof}


We prove the following generalization of Corollary 4.2 in
\cite{B1} and Theorem 4.1 in \cite{CMF}. The techniques are very
similar to the ones used by Friedman and Casalaina-Martin.


\begin{prop}\label{prop:2.4} Let $C$ be a reduced nodal
sextic contained in $\mathbb{P}^2$ and let $\theta$ be a
theta-characteristic on $C$ such that $h^0(C,\theta)=1$. Chosen
homogeneous coordinates $x_1,x_2,x_3$ in $\mathbb{P}^2$, there
exists a matrix
\begin{eqnarray}
M=\left(
\begin{array}{cccc}
l_{11} & l_{12} & l_{13} & q_1 \\
l_{21} & l_{22} & l_{23} & q_2 \\
l_{31} & l_{32} & l_{33} & q_3 \\
q_1 & q_2 & q_3 & f
\end{array}
\right),
\end{eqnarray}
where $l_{ij}$, $q_k$ and $f$ are polynomials respectively of
degree 1, 2 and 3 in $x_1$, $x_2$ and $x_3$ and $l_{ij}=l_{ji}$,
such that the following sequence is exact:
\begin{eqnarray}
0\longrightarrow\mathcal{O}_{{\mathbb{P}}^2}(-2)^3\oplus\mathcal{O}_{{\mathbb{P}}^2}(-3)\stackrel{M}{\longrightarrow}\mathcal{O}_{{\mathbb{P}}^2}(-1)^3\oplus\mathcal{O}_{{\mathbb{P}}^2}\longrightarrow\theta\longrightarrow
0.
\end{eqnarray}
In particular, $\det M$ is the equation which defines $C$.

Conversely, if $C$ is a reduced nodal plane sextic and $M$ is a
matrix of type (1) fitting in a short exact sequence as (2), then
$\theta$ is a theta-characteristic and $h^0(\theta)=1$.
\end{prop}


\begin{proof} Let $\theta=\nu_*L$, where $(Y,L)$ is a stable spin
curve whose stable model is $C$ and let $\nu:Y\rightarrow C$ be
the contraction map. Let $Z:=\overline{Y-(\cup_{i\in I} E_i)}$ and
$\tilde L:=L|_Z$, where $\{E_i:i\in I\}$ is the set of the
irreducible components of $Y$. By Definition \ref{def:2.1}, we
have $(\tilde L)^{\otimes 2}=\omega_Z$.

Now we get the following isomorphisms:
\[
\begin{array}{rl}
\mathcal{H}om(\theta,\omega_C) & \cong\mathcal{H}om(\nu_*\tilde
L,\omega_C)\nonumber \\
 & \cong\nu_*(\tilde L^{-1}\otimes\mathcal{O}(-\nu^{-1}(S_\theta))\otimes\nu^*\omega_C)\nonumber \\
 & \cong\nu_*(\tilde L^{-1}\otimes\omega_{Z})\nonumber \\
 & \cong\nu_*\tilde L\cong\theta,\nonumber
\end{array}
\]
where the second isomorphism follows from the same calculation as
in Theorem 4.1 in \cite{CMF}.

Adapting the proof of Theorem B in \cite{B1} and using the same
calculations of the proof of Corollary 4.2 in \cite{B1} we get the
desired result (the previous isomorphism shows that the matrix $M$
is symmetric). The converse follows using the results in
\cite{B1}.\end{proof}


For a reduced nodal plane sextic $C$ in $\mathbb{P}^5$, let
$\Pi(C)\subset\mathbb{P}^5$ be the plane containing $C$. Let $X$
be a cubic 4-fold in ${\mathbb{P}}^5$ containing a plane $P$ and
let $C$ be a reduced nodal plane sextic in $\mathbb{P}^5$ such
that $P\cap\Pi(C)=\emptyset$. Assume that $\theta$ is a
theta-characteristic on $C$ so that $h^0(C,\theta)=1$.

\begin{definition}\label{def:assdisc} We say that $(X,P)$ is {\it associated to} $(C,\theta)$ or that $(C,\theta)$ is a {\it discriminant curve} of $(X,P)$ if
there exist homogeneous coordinates $x_1,x_2,x_3,u_1,u_2,u_3$ in
${\mathbb{P}}^5$ such that:
\begin{itemize}
\item[(a)] the equations of $P$ are $x_1=x_2=x_3=0$, the equations of $\Pi(C)$ are
$u_1=u_2=u_3=0$ while the equation of $X$ is
\[
F:=\sum_{i,j=1,2,3}l_{ij}u_iu_j+2\sum_{k=1,2,3}q_ku_k+f,
\]
where $l_{ij}, q_k$ and $f$ are polynomials in $x_1$, $x_2$ and
$x_3$ of degree 1, 2 and 3 respectively;
\item[(b)] if $M(X,P,C)$ is the matrix
\[
M(X,P,C)=\left(
\begin{array}{cccc}
l_{11} & l_{12} & l_{13} & q_1 \\
l_{21} & l_{22} & l_{23} & q_2 \\
l_{31} & l_{32} & l_{33} & q_3 \\
q_1 & q_2 & q_3 & f
\end{array}
\right),
\]
whose coefficients are the polynomials $l_{ij}$, $q_k$ and $f$ as
in (a), then $\theta$ and $M(X,P,C)$ fit in a short exact sequence
\[
0\longrightarrow\mathcal{O}_{\Pi(C)}(-2)^3\oplus\mathcal{O}_{\Pi(C)}(-3)\stackrel{M(X,P,C)}{\longrightarrow}\mathcal{O}_{\Pi(C)}(-1)^3\oplus\mathcal{O}_{\Pi(C)}\longrightarrow\theta\longrightarrow
0
\]
of sheaves on $\Pi(C)$ of type (2);
\item[(c)] the equations of $C$ in $\mathbb{P}^5$ are
$\mathrm{det}M(X,P,C)=u_1=u_2=u_3=0$.
\end{itemize}\end{definition}

We will often say that a cubic 4-fold $X$ containing a plane $P$
is {\it associated} to a plane reduced nodal sextic $C$ in
$\mathbb{P}^5$ if there exists a theta-characteristic $\theta$ on
$C$ so that $h^0(C,\theta)=1$ and $(X,P)$ is associated to
$(C,\theta)$. Furthermore, a reduced nodal plane sextic
$C\subset\Pi(C)\subset\mathbb{P}^5$ is a {\it discriminant curve}
of a cubic 4-fold $X$ containing a plane $P$ is there is a
theta-characteristic $\theta$ on $C$ such that $h^0(C,\theta)=1$
and $(C,\theta)$ is a discriminant curve of $(X,P)$.

Let $X\subset\mathbb{P}^5$ be a cubic 4-fold containing a plane
$P$ and let $C$ be a reduced nodal plane sextic in $\mathbb{P}^5$.
Assume that $\theta$ is a theta-characteristic on $C$ such that
$h^0(C,\theta)=1$ and that $(C,\theta)$ is a discriminant curve of
$(X,P)$. We define the rational map
\[
\pi_{P,C}:X\dashrightarrow\Pi(C)
\]
as the projection from the plane $P$ onto $\Pi(C)$. Given a point
$p\in\Pi(C)$ and the 3-dimensional projective space
${\mathbb{P}}^3_p$ containing $p$ and $P$, we define
$F_p:={\mathbb{P}}^3_p\cap X$. We have $F_p=Q_p\cup P$, where
$Q_p$ is a quadric surface which is called the {\it fiber} of the
projection $\pi_{P,C}$ over $p$.

\begin{remark}\label{rmk:interpr} We can give a more geometric
interpretation of a discriminant curve $C$ of a cubic 4-fold $X$
containing a plane $P$. Indeed, it is a classical result (see
\cite{V}) that if $X$ and $C$ are smooth, the curve $C$ can be
thought as the discriminant curve of the quadric bundle
$\pi:\mathrm{Bl}_P(X)\rightarrow\Pi(C)$ obtained from
$\pi_{P,C}:X\dashrightarrow\Pi(C)$ by blowing up $P$ inside $X$.
Proposition \ref{lem:3.4} and Lemmas \ref{lem:3.3} and
\ref{lem:3.5} show that the same holds true when $C$ is a nodal
reduced plane sextic.\end{remark}

We write $W(X,P,C)$ for the net of conics given by
\[
C_p:=Q_p\cap P,
\]
when $p$ varies in $\Pi(C)$. $B(X,P,C)$ is the base points locus
of $W(X,P,C)$. To the pair $(C,\theta)$ it is naturally associated
a cubic plane curve $D\subset\Pi(C)$. If $M$ is the matrix with
polynomial coefficients defined as in Equations (1) and (2), the
equation of $D$ in $\Pi(C)$ is $\det G=0$, where
\[
G:=\left(
\begin{array}{ccc}
l_{11} & l_{12} & l_{13} \\
l_{21} & l_{22} & l_{23} \\
l_{31} & l_{32} & l_{33}
\end{array}
\right)
\]
is $3\times 3$ minor corresponding to the linear part of $M$. We
introduce the sets
\[
I_C:=D\cap C \qquad \qquad S_C:=\Sing(C)-I_C
\]
and
\[
\tilde S_\theta:=\{p\in\Sing(C):p\in \supp
(\theta/s\mathcal{O}_{C})\},
\]
where $s$ is such that $\langle s \rangle =H^0(C,\theta)$.

\begin{remark}\label{rmk:constr} Let $(C,\theta)$ be a discriminant curve of $(X,P)$, where
$X\subset\mathbb{P}^5$ is a cubic 4-fold, $P\subset X$ is a plane,
$C$ is a reduced nodal plane sextic in $\mathbb{P}^5$ and $\theta$
is a theta-characteristic on $C$ such that $h^0(C,\theta)=1$.
Assume that $P_1,P_2\subset X$ are two planes so that $P_1\cap
P_2=\emptyset$. Consider $\varphi\in\mathrm{PGL}(6)$ such that
$\varphi(P_1)=P$ and $\varphi(P_2)=\Pi(C)$. By Definition
\ref{def:assdisc}, $(\varphi^{-1}(C),\varphi^*(\theta))$ is a
discriminant curve of $(X,P_1)$. Moreover, any two discriminant
curves $(C_1,\theta_1)$ and $(C_2,\theta_2)$ of $(X,P)$ are
isomorphic.\end{remark}


\section{Singular points of cubic 4-folds containing a plane}


This section is devoted to prove the following theorem which
describes the singular points of a cubic 4-fold containing a plane
which is associated to a reduced plane sextic with at most nodal
singularities. This results also gives an estimate of the rank of
the free abelian group $\mathrm{NS}_2(X)$ generated by the
equivalence classes of the algebraic cycles of codimension 2 in
$X$. We use the notations introduced in Section
\ref{sec:thetacubic}.


\begin{thm}\label{thm:3.1} Let $X$ be a cubic 4-fold
containing a plane $P$ with discriminant curve $(C,\theta)$, where
$C$ is a nodal reduced plane sextic and $\theta$ is a
theta-characteristic on $C$ such that $h^0(C,\theta)=1$. Then all
the singular points of $X$ are double points, the set $\Sing(X)$
is zero dimensional and
\[
\begin{array}{rl}
\Sing(X) & =\left(\bigcup_{p\in(\Sing(C)-\tilde S_\theta)}\Sing(Q_p)\right)\coprod B(X,P,C)\nonumber \\
 & =\left(\bigcup_{p\in S_C}\Sing(Q_p)\right)\coprod
B(X,P,C),\nonumber
\end{array}
\]
where $Q_p$ is the fiber over $p\in\Pi(C)$ of the projection
$\pi_{P,C}$. In particular, $X$ is irreducible and
\[
\#S_C \leq \#\Sing(X)\leq \#S_C+3.
\]
Moreover, the free abelian group $\mathrm{NS}_2(X)$ generated by
the equivalence classes of the algebraic cycles of codimension 2
in $X$ contains $2(\#S_\theta)+1$ distinct classes represented by
planes and
\[
\mathrm{rk}\mathrm{NS}_2(X)\geq\#S_\theta+2.
\]\end{thm}


The proof of this theorem requires some preliminary results. We
start with the following lemmas and propositions.


\begin{lem}\label{lem:3.2} Let $X$ be a cubic 4-fold
containing a plane $P$ associated to a reduced nodal plane sextic
$C$ with a theta-characteristic $\theta$ such that
$h^0(C,\theta)=1$. Let $x\in\Sing(X)-P$. If
$x_1,x_2,x_3,u_1,u_2,u_3$ are homogeneous coordinates in
$\mathbb{P}^5$ as in Definition \ref{def:assdisc} and if
$x=(a_1:a_2:a_3:b_1:b_2:b_3)$, then
\[
p:=\pi_{P,C}(x)=(a_1:a_2:a_3)\in\Sing(C).
\]
Furthermore, if $f$, $q_1$, $q_2$ and $q_3$ are the polynomials
appearing in the equation of $X$ as in Definition
\ref{def:assdisc}, then
\[
f(p)=q_1(p)=q_2(p)=q_3(p)=0
\]
and all the first partial derivatives of $f$ evaluated in $p$ are
zero.\end{lem}


\begin{proof} For $a_j \neq 0$, the map $g\in$PGL$(6)$ defined by
\[
(x_1:x_2:x_3:u_1:u_2:u_3)\mapsto(x_1:x_2:x_3:u_1-(b_1/a_j)x_j:u_2-(b_2/a_j)x_j:u_3-(b_3/a_j)x_j)
\]
is such that $g(x)=(a_1:a_2:a_3:0:0:0)$. Obviously, $g|_P=
\text{Id}$ and $\pi_{P,C}(x)=\pi_{P,C}(g(x))=p$ while
$P':=g(\Pi(C))$ can be different from $\Pi(C)$. Thus, up to
passing to the discriminant curve $(g(C),g_*(\theta))$ (see Remark
\ref{rmk:constr}), we can suppose without loss of generality that
$x=(a_1:a_2:a_3:0:0:0)$.

Evaluating the equation of $X$ in $x$ and calculating its partial
derivatives in $x$ in these new coordinates, we get the following
equalities:
\[
f(a_1,a_2,a_3)=q_1(a_1,a_2,a_3)=q_2(a_1,a_2,a_3)=q_3(a_1,a_2,a_3)=0,
\]
where $f$, $q_1$, $q_2$ and $q_3$ are the polynomials appearing in
the equation of $X$ as in Definition \ref{def:assdisc}. Moreover,
the point $(a_1:a_2:a_3:0:0:0)$ is a singular point for the cubic
$f=0$ in $\Pi(C)$. Thus all the first partial derivatives of $f$
in $p$ are zero.

If $h$ is the equation of $C$ in the local ring of $\Pi(C)$ at
$p$, we get $h=gf$ modulo $m_p^2$, where $m_p$ is the maximal
ideal and $g$ is a polynomial of degree 3.  The cubic $f=0$ is
singular in $p$ and thus $h=0$ modulo $m_p^2$. This implies that
$p$ is a singular point of $C$ (similar calculations are done in
\cite{CMF} for cubic 3-folds).
\end{proof}


\begin{lem}\label{lem:3.3} Let $X$ be a cubic 4-fold
containing a plane $P$ associated to a reduced nodal plane sextic
$C$ with a theta-characteristic $\theta$ such that
$h^0(C,\theta)=1$.

{\rm (i)} If $p\in\tilde S_\theta-S_\theta$, then $Q_p$ is a cone
(i.e $\Sing(Q_p)$ consists of one point) and $\Sing(Q_p)\subset P$.

{\rm (ii)} If $p\in\Sing(C)-\tilde S_\theta$, then $Q_p$ is a cone
and $Q_p\cap P$ is a smooth conic.\end{lem}


\begin{proof} Let $x_1,x_2,x_3,u_1,u_2,u_3$ be homogeneous coordinates in
$\mathbb{P}^5$ satisfying (a)--(c) in Definition
\ref{def:assdisc}. By Definition \ref{def:2.2}, $\theta=\nu_* L$,
where $L$ is a line bundle defined on a stable spin curve $Y$
whose stable model is $C$ and $\nu$ is the contraction map.
Moreover, $C$ is obtained contracting all the exceptional
components of $Y$ and ${S_\theta}$ is the image in $C$ of all
these components. If $p\not\in S_{\theta}$, $\theta$ is locally
invertible in $p$ (see Lemma \ref{lem:2.1}). This implies that if
$p\in\Sing(C)-\tilde S_\theta$ or $p\in\Sing(C)-S_\theta$, the
rank of $M(X,P,C)(p)$ must be 3 (where $M(X,P,C)(p)$ means the
matrix $M(X,P,C)$ evaluated in $p$).

For any $p=(a:b:c:0:0:0)\in\Pi(C)$, let ${\mathbb{P}}^3_p$ be the
3-dimensional projective space generated by $p$ and $P$. Chosen in
$\mathbb{P}_p^3$ the homogeneous coordinates $u_1,u_2,u_3,t$, the
equation of $F_p=\mathbb{P}_p^3\cap X$ in $\mathbb{P}_p^3$ is
\begin{eqnarray}
L:=\sum_{i,j=1,2,3}l_{ij}(a,b,c)tu_iu_j+2\sum_{k=1,2,3}q_k(a,b,c)t^2u_k+f(a,b,c)t^3=0.
\end{eqnarray}
When $p\in\Sing(C)-S_\theta$, the matrix $M(X,P,C)(p)$ has rank 3.
Hence, by Equation (3), the quadric $Q_p\subset F_p$ is a cone and
it has only one singular point. This proves the first part of (i)
and (ii).

Let us consider the cubic plane curve $D\subset\Pi(C)$ defined in
the last part of Section \ref{sec:thetacubic}. Its equation in
$\Pi(C)$ is the determinant of the matrix $G$ which is the linear
part of $M(X,P,C)$. In \cite{CMF} it is proved (in the case of a
cubic 3-fold, but their proof can be easily modified for our
purposes) that $p\in D\cap C$ if and only if $p$ is in the support
of $\theta/\mathcal{O}_C$ (the quotient is given by the inclusion
$\mathcal{O}_X\hookrightarrow\theta$, defined by a section $s$ of
$\theta$ such that $H^0(C,\theta)=\langle s \rangle$). This means
that if $G(p)\equiv 0$, then $\theta$ is not locally isomorphic to
$\mathcal{O}_X$ in $p$. By definition, this would imply
$p\in\tilde S_\theta$. Thus, if $p\in\Sing(C)-\tilde S_\theta$,
then $G(p)\not\equiv 0$ and, more precisely, rk$G(p)=3$. On the
other hand, if $p\in\tilde S_\theta-S_\theta$, by definition,
$\det G(p)=0$ and thus the rank of $G(p)$ is less than three.

The equation of the conic $Q_p\cap P$ is
\[
\sum_{i,j}l_{ij}(a,b,c)u_iu_j=0.
\]
Hence, if $p\in\Sing(C)-\tilde S_\theta$, then the curve $Q_p\cap
P$ is smooth since the rank of $M(X,P,C)(p)$ and of $G(p)$ is 3.
If $p\in\tilde S_\theta-S_\theta$, then $\Sing(Q_p)\subset P$,
because the rank of $G(p)$ is smaller than 3.\end{proof}


\begin{prop}\label{lem:3.4} Let $X$ be a cubic 4-fold
containing a plane $P$ associated to a reduced nodal plane sextic
$C$ with a theta-characteristic $\theta$ such that
$h^0(C,\theta)=1$. Let $p\in\Sing(C)-\tilde S_\theta$. Then
\[
\Sing(X)\supseteq\Sing(Q_p).
\]
In particular, if $\Sing(C)\neq \tilde S_\theta$, then
$\Sing(X)\neq\emptyset$.
\end{prop}


\begin{proof} As usual, assume that $x_1,x_2,x_3,u_1,u_2,u_3$ are homogeneous coordinates in
$\mathbb{P}^5$ satisfying (a)--(c) in Definition
\ref{def:assdisc}. We prove this result in three steps.

\vspace{0,1cm}

\noindent {\sc Step 1.} If $p\in C-\Sing(C)$ then
rk$M(X,P,C)(p)=3$, where $M(X,P,C)(p)$ means the matrix $M(X,P,C)$
evaluated at $p$. Moreover, $\Sing (Q_p)$ consists of one point
and $Q_p$ is a cone.


\begin{proof} This follows from the trivial remark that if
rk$M(X,P,C)(p)<3$, then the determinants of all the $3\times 3$
minors of $M(X,P,C)$ are zero in $p$. Moreover, the partial
derivatives of the equation of $C$ is $\sum_i d_i m_i$, where
$d_i$ is the derivative of one of the polynomials entries of
$M(X,P,C)$ and $m_i$ is the determinant of a $3\times 3$ minor of
$M(X,P,C)$.\end{proof}


\noindent {\sc Step 2.} If $p\in\Sing(C)-\tilde S_\theta$ then
there is an open neighborhood $U$ of $p$ in $C$ such that the map
$\varphi:U\rightarrow X$ defined by
\[
\varphi:q'\longmapsto \Sing(Q_{q'})
\]
is an isomorphism onto its image and
$\pi_{P,C}|_{\varphi(U)}:\varphi(U)\rightarrow U$ is an
isomorphism.


\begin{proof} We start considering the simplest case of a smooth
plane sextic $C$. It is possible to introduce the embedding
$\psi_{|\theta(1)|}:C\longrightarrow {\mathbb{P}}^5$, defined by
the linear system $|\theta(1)|$. It can be geometrically described
by taking the singular points of the singular quadrics in the
fibration $\pi_{P,C}:X\dashrightarrow\Pi(C)$ (see, for example,
\cite{V} and \cite{B2}). In particular, the image in
${\mathbb{P}}^5$ of a smooth sextic is smooth itself.

We can produce an analogous construction in the reduced nodal
case. Indeed, let $U$ be a suitable open neighborhood of $p$. By
Step 1, if $p\in U-\Sing(C)$, then $Q_p$ has only one singular
point. By Lemma \ref{lem:3.3}(ii) the same is true when
$p\in\Sing(C)-\tilde S_\theta$. Thus, at least, $\varphi$ is
well-defined.

The injectivity follows easily. Indeed, let $q_1, q_2\in U$ be
such that
\[
\{x\}=\Sing(Q_{q_1})=\Sing(Q_{q_2}),
\]
where $Q_{p_1}\not\equiv Q_{p_2}$. Then there is a plane $P_1$
meeting $P$ along a line and such that
\[
C':=Q_{q_1}\cap
P_1 \qquad \text{and} \qquad C'':=Q_{q_2}\cap
P_1
\]
are irreducible plane conics. This means that $F_{q_1}$ would
contain $C'$, $C''$ and $x$. Thus $F_{q_1}\supset Q_{q_1}\cup
Q_{q_2}$, where $Q_{q_1}\not\equiv Q_{q_2}$. This is absurd, since
$\overline{F_{q_1}-P}$ is of degree 2.

Now we can show that $\varphi$, restricted to a suitable open
neighborhood $U$ is an isomorphism onto its image and that
$\pi_{P,C}|_{\varphi(U)}$ maps $\varphi(U)$ isomorphically onto
$U$. By Lemma \ref{lem:3.3}, when $p\in\Sing(C)-\tilde S_\theta$,
there is an open neighborhood $V$ of $p$ such that
$\Sing(Q_q)\not\subseteq P$, when $q\in V$. Now we recall that the
equation of $Q_q$ inside the projective space ${\mathbb{P}}^3_q$
containing $q$ and $P$ is obtained by dividing up Equation (3) by
$t$. Thus, by looking at the partial derivatives of (3), we see
that there is an open neighborhood $W$ of $p\in\Sing(C)-\tilde
S_\theta$ such that the map $\varphi$ is given by
\[
(a,b,c)\longmapsto
(a,b,c,\psi_1(a,b,c),\psi_2(a,b,c),\psi_3(a,b,c)),
\]
where $\psi_1$, $\psi_2$ and $\psi_3$ are rational functions
defined in $W$. Restricting $W$ to $C$ we get the thesis.
\end{proof}


\noindent {\sc Step 3.} If $q\in\Sing(C)-\tilde S_\theta$ then
$\Sing(Q_q)\subseteq\Sing(X)$.


\begin{proof} Let us consider the differential
\[
(\mathrm{d}\pi_{P,C})_{q'}:T_{q'}\varphi(U)\longrightarrow
T_{q}{\mathbb{P}}^2,
\]
where $q':=\varphi(q)$. It has rank 2 in $q'$. Moreover, by
definition, $q'$ is the singular point of the quadric $Q_q$ which
is the fiber of the projection $\pi_{P,C}$ over the node $q$. This
means that $\dim T_{q'} Q_q=3$. Since
\[
\dim T_{q'}X=\dim T_{q'}Q_q+\dim (\mathrm{Im}(d\pi_{q'}))=3+2=5,
\]
$q'$ is singular in $X$ and this proves Step 3.\end{proof}


This concludes the proof of the lemma. \end{proof}


\begin{lem}\label{lem:3.5} Assume that $X$ is a cubic 4-fold
containing a plane $P$ associated to the pair $(C,\theta)$, where
$C$ is a reduced nodal plane sextic, $\theta$ is a
theta-characteristic such that $h^0(C,\theta)=1$ and $p\in
S_\theta$. Then

{\rm (i)} $\Sing(X)\cap Q_p\subset P$;

{\rm (ii)} $Q_p=P_1\cup P_2$, where $P_1$ and $P_2$ are distinct
planes.

\noindent Moreover, $S_\theta\subseteq\tilde S_\theta$.
\end{lem}


\begin{proof} Looking at the short exact
sequence of sheaves in $\Pi(C)\cong\mathbb{P}^2$
\begin{eqnarray}
0\longrightarrow\mathcal{O}_{\Pi(C)}(-2)^3\oplus\mathcal{O}_{\Pi(C)}(-3)\stackrel{M(X,P,C)}{\longrightarrow}\mathcal{O}_{\Pi(C)}(-1)^3\oplus\mathcal{O}_{\Pi(C)}\longrightarrow\theta\longrightarrow
0,
\end{eqnarray}
and at Definition \ref{def:2.2}, we see that if $p\in S_\theta$,
then $\mathrm{rk}M(X,P,C)(p)=2$. In particular, $Q_p$ is the union
of two distinct planes as stated in (ii).

Suppose that $x_1,x_2,x_3,u_1,u_2,u_3$ are homogeneous coordinates
in $\mathbb{P}^5$ satisfying (a)--(c) in Definition
\ref{def:assdisc}. Let $x\in\Sing(X)$ be such that
$p:=\pi_{P,C}(x)\in S_\theta$. By Lemma \ref{lem:3.2}, if
$x\in\Sing(X)$ then $p=\pi_{P,C}(x)\in\Sing(C)$ and the
polynomials $f$, $q_i$ (for $i=1,2,3$) and all the first partial
derivatives of $f$ are zero in $p$ (see Definition
\ref{def:assdisc} and Lemma \ref{lem:3.2} for the definitions of
the polynomials $f$ and $q_i$). If $x\not\in P$ but $p\in
S_\theta$, the previous remarks implies that, up to a change of
variables, we can suppose that the matrix $M(X,P,C)(p)$ (which has
rank 2) is
\[
\left(
\begin{array}{cccc}
1 & 0 & 0 & 0 \\
0 & 1 & 0 & 0 \\
0 & 0 & 0 & 0 \\
0 & 0 & 0 & 0
\end{array}
\right).
\]
In particular, the equation of the curve $C$, in the local ring
$R$ of $\Pi(C)$ at $p$, would be congruent to $-l_{12}^2$ mod
$m_p^3$, where $m_p$ is the maximal ideal of $R$. This gives the
desired contradiction since $C$ has no cusp of multiplicity 2
(similar calculations are done in \cite{CC} for cubic 3-folds in a
slightly different context). Hence $x\in P$ and (i) holds.

The last statement in the lemma is a consequence of the fact that,
since $\theta$ is the push-forward of a line bundle on a stable
spin curve (see Definition \ref{def:2.1}), the sheaf $\theta$ is
not locally isomorphic to $\mathcal{O}_C$ in $p$ (by Lemma
\ref{lem:2.1}), for every $p\in S_\theta$. Thus $p$ is in the
support of $\theta/s\mathcal{O}_C$, where $s$ is such that
$\langle s \rangle=H^0(C,\theta)$. This implies that $p$ is in
$\tilde S_\theta$. \end{proof}


\begin{remark}\label{rmk:3.1} We can find examples which
show that $S_\theta$ can be different from $\tilde S_\theta$. One
can consider the case of an irreducible plane sextic $C$ with a
node $p$ and a theta-characteristic $\theta$ on $C$ such that
$H^0(C,\theta)=\langle s\rangle$, the section $s$ is a cubic
intersecting $C$ in $p$ and $\theta$ does not come from a
theta-characteristic on the total normalization of $C$. It is easy
to see that, in this case, $S_\theta=\emptyset$ and
$\tilde{S}_\theta=\{p\}$.

We can describe some more interesting examples. Indeed, let
$x_1,x_2,x_3,u_1,u_2,u_3$ be homogeneous coordinates in
$\mathbb{P}^5$. Assume that $C$ is the union of the nodal plane
cubic $C_1$ whose equations are
$u_1=u_2=u_3=x_2^2x_3-x_1^3+x_1^2x_3=0$ and of a smooth plane
cubic $C_2$ whose equations are $u_1=u_2=u_3=f=0$. Moreover,
suppose that $C_1\cap C_2$ consists of nine points distinct from
the nodes of $C_1$. There exists a theta-characteristic $\theta_1$
on $C_1$ such that $h^0(C_1,\theta_1)=0$. In particular, we can
write the equation of $C_1$ in $\Pi(C)$ as the determinant of the
matrix
\[
M_1:=\left(
\begin{array}{ccc}
0 & x_1 & x_2 \\
x_1 & -x_3 & 0 \\
x_2 & 0 & x_1 + x_3
\end{array}
\right).
\]
The equation of $C$ in $\Pi(C)$ is then the determinant of the
matrix
\[
M_2:=\left(
\begin{array}{cc}
M_1 & 0 \\
0 & f
\end{array}
\right).
\]
As we have seen in the proof of Lemma \ref{lem:3.2},
$p\in\Sing(C)$ is in $\tilde S_\theta$ if and only if $p\in C\cap
D$. This means that the point $p:=(0:0:1:0:0:0)$ is in $\tilde
S_\theta$ but not in $S_\theta$. This is clear since the rank of
$M_2$ in $p$ is 3 (simply because $C_2$ does not intersect $C_1$
in its node) but the intersection of the quadric $Q_p$ with the
plane $P$ is the union of two lines (because the rank of $M_1$ in
$p$ is 2).\end{remark}




\begin{definition}\label{def:3.1} A \emph{couple of planes} is given
by two planes $P_1$ and $P_2$ such that
\[
P_1\cap P_2=r,
\]
where $r$ is a line. We write $(P_1,P_2)$ for such a couple.
\end{definition}


We can prove the following proposition.


\begin{prop}\label{prop:3.6} Let $X$ be a cubic 4-fold
containing a plane $P$ associated to a discriminant curve
$(C,\theta)$. If $q_1,q_2\in S_\theta$ and $Q_{q_1}=P_1\cup P_2$,
$Q_{q_2}=S_1\cup S_2$, where $P_i$ and $S_j$ are planes, then
\[
P_i\cap S_j=1\mbox{ point},
\]
for any $i,j\in\{1,2\}$. Furthermore, given a couple of planes
$(P_1,P_2)$ in $X$, there exists a plane $P\subset X$ such that
$P_1\cup P_2$ is a fiber of the projection $\pi_{P,C}$.\end{prop}


\begin{proof} Observe that, due to Lemma \ref{lem:3.5}, if $p\in S_\theta$ then $Q_p=P_1\cup
P_2$, where $P_1$ and $P_2$ are planes. Moreover $P_1 \neq P$ and
$P_2 \neq P$. Indeed, assume that $x_1,x_2,x_3,u_1,u_2,u_3$ are
homogeneous coordinates in $\mathbb{P}^5$ as in Definition
\ref{def:assdisc}. Let ${\mathbb{P}}^3_p$ be the projective space of
dimension 3 containing $p=(a:b:c:0:0:0)$ and $P$. As we noticed in
the proof of Lemma \ref{lem:3.3}, the equation of $F_p$ in
${\mathbb{P}}^3_p$ is the polynomial $L$ defined by Equation (3) in
the homogeneous coordinates $u_1,u_2,u_3,t$. The choice of these
coordinates in Lemma \ref{lem:3.3} was such that the equation of $P$
in ${\mathbb{P}}^3_p$ is $t=0$. Thus $P_i\equiv P$ if and only if,
in ${\mathbb{P}}^3_p$, $L=t^2M$, where $M$ is a polynomial of degree
1 in $u_1,u_2,u_3$ and $t$. This happens if and only if the matrix
$G$ (see the end of Section \ref{sec:thetacubic} and Lemma
\ref{lem:3.3} for its definition) is the zero matrix in $p$ (i.e.
$G(p)\equiv 0$). This would imply $\mathrm{rk}M(X,P,C)(p)=1$ but
this is a contradiction since $\mathrm{rk}M(X,P,C)(p)=2$, when $p\in
S_\theta$.

Let $(P_1,P_2)$ and $(S_1,S_2)$ be two couples of planes such that
$Q_{q_1}=P_1\cup P_2$ and $Q_{q_2}=S_1\cup S_2$, for $q_1,q_2\in
S_\theta$. By Lemma \ref{lem:3.5} the intersection of $P_1\cup
P_2$ (and $S_1\cup S_2$) with the plane of projection $P$ is the
union of two lines (not necessarily distinct). Hence $S_i$ and
$P_j$ meet each other at least in one point $P_i\cap S_j\cap P$.
Suppose that there exist $i,j\in\{1,2\}$ such that $S_i\cap P_j$
is a line $l_1$ and let $l_2$ be the line $P_1\cap P_2$.

If $l_1,l_2\subset P$, then $l_1\equiv l_2$ because, otherwise,
$P_j\equiv P$ (recall that $l_1,l_2\subset P_j$). Hence, we would
have a plane $R$ meeting $P$, $S_i$, $P_1$ and $P_2$ along 4 lines
$r_1$, $r_2$, $r_3$ and $r_4$. Thus $F_{q_1}\supset P_1\cup P_2$
would contain $l_1$, $r_1$, $r_2$, $r_3$, $r_4$ and so it would
also contain $S_i(\not\equiv P)$. If $l_1\subset P$ but
$l_2\not\subset P$, then $S_i\cup P_1\cup P_2\subset F_{q_1}$
because $S_i\cup P_j\subset\mathbb{P}^3_{q_1}$. On the other hand,
if $l_1\not\subset P$, then $(S_i\cap P)\cup l_1\subset F_{q_1}$
and $S_i\subset F_{q_1}$. In all these cases we get a
contradiction since $\overline{F_{q_1}-P}$ is a quadric. Hence we
can conclude that, for any $i,j\in\{1,2\}$, $S_i\cap P_j$ consists
of one point.

Suppose now that $S_1$ and $S_2$ are a couple of planes. Let
$P'\subset\mathbb{P}^5$ be a plane such that $P'\cap
S_1=\emptyset$. Consider the projection
$\pi_{S_1,P'}:X\dashrightarrow P'$ from $S_1$. Remark
\ref{rmk:constr} shows that there is a nodal plane sextic $C'$
with a theta-characteristic $\theta'$ such that $(C',\theta')$ is
a discriminant curve of $(X,S_1)$, $P'\equiv\Pi(C')$ and
$\pi_{S_1,P'}$ coincides with the projection
$\pi_{S_1,C'}:X\dashrightarrow\Pi(C')$. The plane $S_2$ is
contained in a fiber of $\pi_{S_1,P'}$. Indeed, since $S_1$ and
$S_2$ meet along a line, there exists a ${\mathbb{P}}^3$
containing both planes and meeting $P'$ in one point. This implies
that $S_2$ is a component of a fiber of
$\pi_{S_1,P'}=\pi_{S_1,C'}$. Applying Lemma \ref{lem:3.5} to the
pairs $(X,S_1)$ and $(C',\theta')$, we see that such a fiber of
$\pi_{S_1,C'}$ must contain two distinct planes. Hence there
exists a plane $S_3$ meeting $S_1$ and $S_2$ along two lines.

This concludes the proof. \end{proof}


\begin{remark}\label{rmk:3.3} Looking more carefully at the proof
of Proposition \ref{prop:3.6} and using Remark \ref{rmk:constr},
we can say that each couple of planes $(P_1,P_2)$ contained in a
cubic 4-fold can actually be described as a triple of planes
$(P_1,P_2,P_3)$ such that $P_i\cap P_j$ is a line for
$i,j\in\{1,2,3\}$ and $i\neq j$.
\end{remark}

\begin{lem}\label{lem:4.3} Let $X$ be a smooth cubic 4-fold
containing a plane $P$ with discriminant curve $(C,\theta)$, where
$C$ is a nodal reduced plane sextic and $\theta$ is a
theta-characteristic on $C$ such that $h^0(C,\theta)=1$. Then
\[
\mathrm{rk}\mathrm{NS}_2(X)\geq\#S_\theta+2.
\]
Moreover, $\mathrm{NS}_2(X)$ contains $2(\#S_\theta)+1$ distinct
classes represented by planes.\end{lem}


\begin{proof} By Proposition
\ref{prop:3.6}, we know that $X$ contains $M:=\#S_\theta$ couples
of planes
\[
(P_{1,1},P_{1,2}),\ldots,(P_{M,1},P_{M,2}).
\]
Let $h^2$, $P$ and $Q$ be the classes in $\mathrm{NS}_2(X)$
corresponding, respectively, to the hyperplane section, the plane
$P$ and the general quadric in the fiber of $\pi_{P,C}$. We recall
that, in $\mathrm{NS}_2(X)$, $ h^2=Q+P$.

If ``$\cdot$'' is the cup product, we have $Q\cdot P=Q\cdot
(h^2-Q)=2-4=-2$ and
\[
P_{i,j}\cdot P=\frac{1}{2}(Q^2\cdot P)=-1,
\]
for $i\in\{1,\ldots,12\}$ and $j\in\{1,2\}$. Hence, since $P\cdot
P=3$, the class of $P_{i,j}$ is distinct from the one of $P$, for
any $i\in\{1,\ldots,M\}$ and $j\in\{1,2\}$. Moreover,
\[
P_{i,1}\cdot
P_{i,2}=\frac{1}{2}\left(Q^2-P_{i,1}^2-P_{i,2}^2\right)=-1
\]
and the classes of $P_{i,1}$ and $P_{i,2}$ are distinct in
$\mathrm{NS}_2(X)$, because $P_{i,j}\cdot P_{i,j}=3$. On the other
hand, by Proposition \ref{prop:3.6}, $P_{i,k}\cdot P_{j,h}=1$, for
$i$ distinct from $j$ and $k,h\in\{1,2\}$. Thus the classes of
$P_{i,k}$ and of $P_{j,h}$ are not the same in $\mathrm{NS}_2(X)$
when $k\neq h$ and $i,j\in\{1,\ldots,M\}$. Summarizing, we have
proved that $P$, $P_{1,1},\ldots,P_{M,1}$ and
$P_{1,2},\ldots,P_{M,2}$ are $2(\#S_\theta)+1$ distinct classes in
$\mathrm{NS}_2(X)$.

Let us consider the sublattice $N\subseteq\mathrm{NS}_2(X)$
defined by:
\[
N:=\langle P,P_{1,2},P_{i,1}:i\in\{1,\ldots,M\}\rangle.
\]
Clearly ${\rm rk}\mathrm{NS}_2(X)\geq{\rm rk}N$. Consider the
matrix $A=(a_{ij})$ representing the intersection form in $N$.
Using the previous calculations, we see that the coefficients of
$A$ are defined as follows (keeping the planes ordered as in the
definition of $N$):
\[
a_{ij}=\left\{\begin{array}{ll} 3 & \mbox{if }i=j;\\ -1 &
\mbox{if }i=1,j>1\mbox{ or }j=1,i>1\mbox{ or }i=2,j=3\mbox{ or }i=3,j=2;\\
1 & \mbox{otherwise.}\end{array}\right.
\]
By an easy calculation we see that the determinant of $A$ is
different from 0. Hence $A$ has maximal rank and the
$\#S_\theta+2$ planes generating $N$ are linearly
independent.\end{proof}

\begin{remark} The estimate in Lemma \ref{lem:4.3} does not depend
on the choice of the discriminant curve $(C,\theta)$. Indeed, by
Remark \ref{rmk:constr}, two discriminant curves $(C,\theta)$ and
$(C',\theta')$ of $(X,P)$ are isomorphic.\end{remark}


We are now ready to prove Theorem \ref{thm:3.1}.

\vspace{0,2cm}

\noindent{\it Proof of Theorem \ref{thm:3.1}.}  For simplicity, we
assume from the beginning homogeneous coordinates
$x_1,x_2,x_3,u_1,u_2,u_3$ in $\mathbb{P}^5$ satisfying (a)--(c) in
Definition \ref{def:assdisc}.

Let $q:=(a:b:c:0:0:0)\in\Pi(C)$ and let ${\mathbb{P}}^3_q$ be the
projective space of dimension 3 containing $q$ and the plane $P$.
Depending on the fact that $a\neq 0$ or $a=0$, such a
3-dimensional space can be described by the equations
\[
bx_1-ax_2=cx_1-ax_3=0 \qquad \text{or} \qquad x_1=cx_2-bx_3=0.
\]
In these two cases, the Jacobian matrix $L_q$ of the equations
describing $F_q=Q_p\cup P$ (see Equation (3)) is either
\begin{eqnarray}
\left(\begin{array}{cccccc} b & -a & 0 & 0 & 0 & 0 \\ c & 0 & -a & 0 & 0 & 0 \\
F_{x_1} & F_{x_2} & F_{x_3} & F_{u_1} & F_{u_2} & F_{u_3}
\end{array}\right) & \text {or} & \left(\begin{array}{cccccc} 1 & 0 & 0 & 0 & 0 & 0 \\ 0 & c & -b & 0 & 0 & 0 \\
F_{x_1} & F_{x_2} & F_{x_3} & F_{u_1} & F_{u_2} & F_{u_3}
\end{array}\right),
\end{eqnarray}
where $F$ is the equation of $X$ as in item (a) of Definition
\ref{def:assdisc} while $F_{x_i}$ and $F_{u_j}$ are the partial
derivatives of $F$ with respect to the variables $x_i$ and $u_j$
($i,j\in\{1,2,3\}$).

From (5) it easily follows that the singular points of $X$ are
contained in the union of the singular points of $F_q$ when $q$
varies in $\Pi(C)$. Since $\Sing(F_q)=\Sing(Q_q)\cup(P\cap Q_q)$,
\[
\Sing(X)\subseteq\left(\bigcup_{p\in
C}\Sing(Q_p)\right)\cup\left(\bigcup_{p\in\Pi(C)}(Q_p\cap
P)\right).
\]

By Lemma \ref{lem:3.2}, if $x\in\Sing(X)-P$ then
$\pi_{P,C}(x)\in\Sing(C)$. Moreover, by Lemma \ref{lem:3.3}(i) and
Lemma \ref{lem:3.5}(i), if $x\in\Sing(X)-P$ then
$\pi_{P,C}(x)\in\Sing(C)-\tilde S_\theta$. Furthermore,
Proposition \ref{lem:3.4} shows that if $p\in\Sing(C)-\tilde
S_\theta$ then $\Sing(Q_p)\subseteq\Sing(X)$. Hence
\[
\Sing(X)\subseteq\left(\bigcup_{p\in(\Sing(C)-\tilde
S_\theta)}\Sing(Q_p)\right)\cup\left(\bigcup_{p\in\Pi(C)}(Q_p\cap
P)\right).
\]

Let $x\in\Sing(X)\cap P$. Then the matrix $L_q(x)$ has rank 2, for
every $q\in\Pi(C)$. This means that $x$ is contained in $P\cap
Q_q$, for every $q\in\Pi(C)$. Thus $x\in B(X,P,C)$. The converse
is also true. Indeed, let $y\in B(X,P,C)$. It is easy to see that
if $e_1=(1:0:0)$, $e_2=(0:1:0)$ and $e_3=(0:0:1)$ in $P$, then
$F_{x_i}(y)=G_{i}(y)$, for $i\in\{1,2,3\}$, where $G_{i}$ is the
equation of the conic $C_{i}=Q_{e_i}\cap P$. Since $y\in
B(X,P,C)\subset P$, $F_{x_i}(y)=F_{u_j}(y)=0$, for
$i,j\in\{1,2,3\}$. Hence $y\in\Sing(X)$ and
\[
\Sing(X)=\left(\bigcup_{p\in(\Sing(C)-\tilde
S_\theta)}\Sing(Q_p)\right)\cup B(X,P,C).
\]
We showed in the proof of Lemma \ref{lem:3.3} that
$S_C=\Sing(C)-\tilde S_\theta$ while Lemma \ref{lem:3.3} proves
that $\Sing(Q_q)\not\subseteq P$, when $q\in\Sing(C)-\tilde
S_\theta$. This implies that
\[
\Sing(X)=\left(\bigcup_{p\in S_C}\Sing(Q_p)\right)\coprod
B(X,P,C).
\]

If $o\in X$ is a singular point with multiplicity $\geq 3$, then
$o\in P$. Indeed, if $o\not\in P$, then we can suppose, without loss
of generality, that $o=(a:b:c:0:0:0)$ (as in the proof of Lemma
\ref{lem:3.2}). By simple calculations, for $i,j\in\{1,2,3\}$,
$F_{u_iu_j}(x)=l_{ij}(x)=0$ where $F_{u_iu_j}(x)$ is the second
partial derivative of $F$ with respect to $u_i$ and $u_j$ and
$l_{ij}(x)$ is the $(i,j)$-linear coefficient in $M(X,P,C)$
evaluated in $x$. Moreover, if $p:=\pi_{P,C}(o)$, then
$f(p)=q_1(p)=q_2(p)=q_3(p)=0$ (Lemma \ref{lem:3.2}) and
$M(X,P,C)(p)$ would have rank 0. This is absurd. If $o\in P$, then,
for any $p\in\Pi(C)$, $Q_p$ is singular in $o$ because $F_p$ has
multiplicity greater or equal to $3$ in $o$. In particular, we would
get the absurd conclusion that $\mbox{det}M_X\equiv 0$. Therefore,
each singular point of $X$ has multiplicity at most 2.

In Proposition \ref{prop:3.6} we proved that $P\not\subset Q_p$,
for all $p\in\Pi(C)$. Moreover, the conics $C_q:=Q_q\cap P$ in
$W(X,P,C)$ can not coincide when $q$ varies in $\Pi(C)$ because,
otherwise, $B(X,P,C)=C_q$ and $P-C_q\subseteq\Sing(X)$. This would
contradict the previous description of the singular locus of $X$.
Hence $B(X,P,C)$ is not a conic.

If $B(X,P,C)$ contains a line $l$ then all the conics in
$W(X,P,C)$ are reducible. On the other hand, the smooth part of
$C$ contains at least one point $q$ such that $C_q=Q_q\cap P$ is a
smooth conic (this easily follows, as in the smooth case, since
the rank of the matrix $M(X,P,C)$ calculated in this point is 3.
Lemma \ref{lem:3.3} gives directly this result for the points in
$\tilde S_\theta$). Hence $B(X,P,C)$ contains a finite number of
points. Obviously, such a number is at most 4.

Let us suppose
\[
B(X,P,C)=\{p_1,p_2,p_3,p_4\}.
\]
Given a triple $\{p_{i_1},p_{i_2},p_{i_3}\}$, they are in general
position. Indeed, if $l= \langle p_{i_1},p_{i_2},p_{i_3} \rangle $
is a line, then each conic in $W(X,P,C)$ would pass through
$p_{i_1}$, $p_{i_2}$ and $p_{i_3}$ and it would contain $l$. Then
$l\subseteq B(X,P,C)$ which contradicts the fact that $B(X,P,C)$
is zero-dimensional.

Let us define the following six lines in $P$:
\[
r_{ij}:= \langle p_i,p_j \rangle,
\]
where $i,j\in\{1,2,3,4\}$, $i\neq j$. The geometric picture is
summarized by the following diagram:

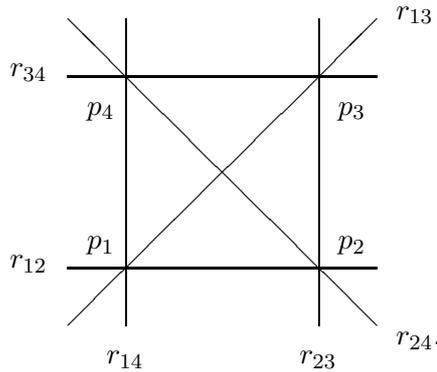
\vspace{3,3cm}

\[
\hspace{-1.9cm}\centerline{\begin{picture}(0,0)\setlength{\unitlength}{1in}
\put(0,0){\line(0,2){1.6}}\put(-0.3,0.3){\line(2,0){1.6}}\put(1,0){\line(0,2){1.6}}\put(-0.3,1.3){\line(2,0){1.6}}
\put(-0.3,0){\line(1,1){1.6}}\put(1.3,0){\line(-1,1){1.6}}
\put(-0.2,0.4){$p_1$}\put(1.1,0.4){$p_2$}\put(-0.2,1.1){$p_4$}\put(1.1,1.1){$p_3$}
\put(-0.6,0.3){$r_{12}$}\put(-0.1,-0.2){$r_{14}$}\put(0.9,-0.2){$r_{23}$}\put(-0.6,1.3){$r_{34}$}
\put(1.4,-0.1){$r_{24}.$}\put(1.4,1.6){$r_{13}$}
\end{picture}}
\]

\vspace{0,3cm}

Let $D\subset\Pi(C)$ be the plane cubic defined in Section
\ref{sec:thetacubic}. The multiplicity of intersection of the
points in the set $D\cap C$ can be at most 2 (see \cite{CMF}).
Hence $D\cap C$ contains at least nine points and $C_p:=Q_p\cap P$
is reducible for $p\in D\cap C$.

If $q_1,q_2\in D\cap S_\theta$, Proposition \ref{prop:3.6}(i)
shows that $C_{q_1}$ and $C_{q_2}$ are distinct because
$Q_{q_1}\cap Q_{q_2}$ consists of a finite number of points. If
$q_1,q_2\in D\cap C$ but $q_1,q_2\not\in S_\theta$, then, by the
same argument as in the proof of Step 2 of Proposition
\ref{lem:3.4} (see the proof of the injectivity of the map
$\varphi$),
\[
\Sing(Q_{q_1})\neq\Sing(Q_{q_2}).
\]
Since $\Sing(Q_{q_i})=\Sing(C_{q_i})$, $C_{q_1}$ and $C_{q_2}$ are
distinct. Thus $W(X,P,C)$ must contain at least five reducible
distinct conics. This gives a contradiction because $W(X,P,C)$
contains only three reducible conics, given by:
\begin{align*}
V_1 & := r_{13}\cup r_{24};\\
V_2 & := r_{12}\cup r_{34};\\
V_3 & := r_{14}\cup r_{23}.
\end{align*}
Hence $B(X,P,C)$ contains at most 3 points and the inequalities
about the number of points in $\Sing(X)$ stated in Theorem
\ref{thm:3.1} hold true.

Furthermore, since $B(X,P,C)$ is finite, $\Sing(X)$ is
zero-dimensional and $X$ is irreducible because the set $\Sing(X)$
is finite. The last part of Theorem \ref{thm:3.1} is exactly Lemma
\ref{lem:4.3}.\hfill$\Box$

\vspace{0,2cm}

We will discuss later (see Examples \ref{ex:4.2}(i) and
\ref{ex:4.2}(ii)) the cases of cubic 4-folds whose singularities
coincide with the singular points of $Q_q$ for
$q\in\Sing(C)-\tilde S_\theta$ or with extra points in $P$.




\section{Smoothness and rationality}


It is clear that, given a nodal plane sextic there are many
associated cubic 4-folds. The first result in this paragraph
describes some sufficient conditions on the geometry of a plane
sextic $C$ such that all the associated cubic 4-folds are
singular. We also give some explicit examples which clarify the
numerical bounds given in Theorem \ref{thm:3.1}.

In the second part of this section, we consider the case of smooth
cubic 4-folds and we introduce a family of smooth rational cubic
4-folds whose discriminant curve is reduced but reducible.

Roughly speaking, Theorem \ref{thm:3.1} proves that the number of
singular points of a cubic 4-folds $X$ containing a plane $P$
depends on the number of nodes of a discriminant curve
$(C,\theta)$ of $(X,P)$ which are not in $S_\theta$. Thus a cubic
4-folds $X$ containing a plane $P$ such that each node of a
discriminant curve $C$ gives rise to a singular point of $X$
realizes the maximal number of singular points compatible with the
geometry of the curve $C$. This is summarized by the following
definition:


\begin{definition}\label{def:4.1} Let $X$ be an irreducible cubic 4-fold containing a plane
$P$. We say that $X$ \emph{realizes the maximal number of singular
points} if any discriminant curve $(C,\theta)$ satisfies the
following conditions:

(i) $C$ is a reduced nodal plane sextic and $\#\Sing(C)>0$;

(ii) $h^0(C,\theta)=1$ and $\# S_\theta=0$.\end{definition}


We say that a reducible plane curve $C$ {\it has general
irreducible components} if, for each irreducible component $C'$
and for each theta-characteristic $\theta$ over $C'$, if $\theta$
is odd (resp. even) then $h^0(C,\theta)=1$ (resp.
$h^0(C,\theta)=0$).


\begin{prop}\label{thm:4.1} {\rm (i)} Let $C$ be a reduced nodal plane sextic in $\mathbb{P}^5$ with general irreducible components and such that either
\begin{itemize}
    \item[(a)] $C$ is irreducible with 10 nodes or

    \item[(b)] $C$ is reducible and it does not contain neither a smooth
cubic nor a quartic with $n$ nodes nor a quintic with $m$ nodes,
where $0\leq n\leq 2$ and $0\leq m\leq 5$.
\end{itemize}
Then there exists at least one singular cubic 4-fold $X$
containing a plane $P$ whose discriminant curve is $C$. Moreover,
all the cubic 4-folds containing a plane associated to $C$ are
singular.

{\rm (ii)} Let $X$ be a cubic 4-fold containing a plane $P$
realizing the maximal number of singular points. Let $C$ be any
discriminant curve of $(X,P)$ with general irreducible components.
Then either $C$ is nodal and irreducible or $C$ is the union of
three smooth conics or of a smooth conic and a quartic (possibly
nodal).\end{prop}


\begin{proof} Let $\theta$ be a theta-characteristic on $C$ such that $h^0(C,\theta)=1$ and let $X$ be a cubic 4-fold containing a plane $P$ associated to $(C,\theta)$.
If $p\in \tilde S_\theta-S_\theta$, using Lemma \ref{lem:3.3}(i) and
the same techniques as in the proof of Steps 2 and 3 of Proposition
\ref{lem:3.4}, one can prove that $\Sing(Q_p)\subseteq\Sing(X)$. In
particular, by Theorem \ref{thm:3.1}, $\Sing(Q_p)\subset B(X,P,C)$.
Analogously, if $x\in\Sing(X)\cap B(X,P,X)$, then
$\pi_{P,C}(x)\in\tilde S_\theta-S_\theta$. Hence all the cubic
4-folds containing a plane associated to a curve $C$ satisfying the
hypotheses of item (i) are singular if $C$ does not have a
theta-characteristic $\theta$ such that
\begin{eqnarray}
\#\Sing(C)=\#
S_\theta\;\;\;\;\;\mbox{and}\;\;\;\;\;h^0(C,\theta)=1.
\end{eqnarray}

Consider the following fact (see \cite{C} or \cite{CC} for the
proof and \cite{CC} for the definitions and the techniques
involving the dual graphs):

\vspace{0,2cm}

\noindent {\sc Fact.} Let $\tilde{C}$ be a stable spin curve whose
stable model is $C$. Let $Z_{\tilde C}:=\overline{\tilde
C-(\cup_{i\in I} E_i)}$, where $\{E_i:i\in I\}$ is the set of the
irreducible components in $\tilde C$, and let $\Gamma_{Z_{\tilde
C}}$ be the dual graph of the curve $Z_{\tilde C}$ (the vertices
of $\Gamma_{Z_{\tilde C}}$ are the irreducible components of
$Z_{\tilde C}$ and its edges are the nodes of $\tilde C$). If
$\Gamma_{Z_{\tilde C}}$ is even (i.e. to each vertex converges an
even number of edges) then the number $M$ of distinct
theta-characteristics over $C$ is such that $M\geq
2^{b_1(\Gamma_{Z_{\tilde C}})}$, where $b_1(\Gamma_{Z_{\tilde
C}})$ is the first Betti number of $\Gamma_{Z_{\tilde C}}$.

\vspace{0,2cm}

\noindent From this and from the easy remark that the push-forward
via the contraction map preserves the parity of the
theta-characteristic, we get two very easy consequences (see
Section 4 in \cite{C} for a brief discussion):
\begin{itemize}
\item[(a.1)] if $2^{b_1(\Gamma_{Z_{\tilde C}})}>1$ then $C$ has an
odd theta-characteristic coming from a line bundle on $\tilde C$
as in Definition \ref{def:2.2};
\item[(b.1)] if $\Gamma_{Z_{\tilde C}}$ has $m$ disjoint components
$Y_1,\ldots,Y_m$, an odd theta-characteristic for $C$ coming from
$\tilde C$ is given by the choice of an odd theta-characteristic
for an odd number of curves $C_i$ corresponding to the graphs
$Y_i$.
\end{itemize}

If $\theta$ satisfies (6), then, by Definitions \ref{def:2.1} and
\ref{def:2.2}, there exists a stable spin curve $(\tilde C,L)$
whose stable model is $C$ and such that $\Sing(C)$ is the image
via the contraction map $\nu:\tilde C\rightarrow C$ of all the
exceptional components of $\tilde C$ and $\theta=\nu_*L$. In this
case, by (b.1), there would be an irreducible component $C_1$ of
$\tilde C$ such that $h^0(C_1,L|_{C_1})=1$. This is impossible,
since when $C$ satisfies the hypotheses of Proposition
\ref{thm:4.1}(i), the irreducible components of such a curve
$\tilde C$ have genus 0.

By the remarks at the beginning of this proof and by (b.1), a
plane sextic $C$ as in item (i) is the discriminant curve of at
least one (singular) cubic 4-fold $X$ containing a plane $P$ if
there are a set $S\subset\Sing(C)$ with $S\neq\Sing(C)$ and a
stable spin curve $(\tilde C,L)$ such that:
\begin{itemize}
\item[(a.2)] $\nu(\tilde{C})=C$;
\item[(b.2)] $S$ is the image via $\nu$ of all the exceptional components of $\tilde C$;
\item[(c.2)] the curve $Z_{\tilde C}:=\overline{\tilde C-(\cup_{i\in I}
E_i)}$ has an irreducible component of arithmetic genus one, where
$\{E_i:i\in I\}$ is the set of the irreducible components in
$\tilde C$.
\end{itemize}
Furthermore, a curve $C$ as in item (i) can only be:
\begin{itemize}
\item[(a.3)] the union of two singular cubics $C_1$ and $C_2$;
\item[(b.3)] the union of a line and of a quintic with 6 nodes;
\item[(c.3)] the union of two lines and of a quartic with 3 nodes;
\item[(d.3)] an irreducible sextic with 10 nodes.
\end{itemize}
In the first case we put $S=\Sing(C_1)\cup(C_1\cap C_2)$. In case
(b.3) we define $S$ as the union of 5 of the 6 nodes of the
quintic and of the 5 intersection points of the line and the
quintic. In the third case, $S$ is the set of all the intersection
points and of 2 of the 3 nodes. In case (d.3), we observe that the
curve $C$ blown up in 9 of the 10 nodes has genus 1. In all of
these cases, due to (a.1)--(b.1), a curve $\tilde C$ satisfying
(a.2)--(c.2) always exists.

For item (ii), observe that if $(X,P)$ realizes that maximal
number of singular points then there exists a discriminant curve
$(C,\theta)$ such that $\theta^2=\omega_C$. Consider the following
cases:

\vspace{0,1cm}

\noindent {\sc Case 1:} $C$ is irreducible with $m$ nodes
($m<11$);

\vspace{0,1cm}

\noindent {\sc Case 2:} $C$ is reducible and it is either the
union of 3 smooth conics or the union of a smooth conic and of a
quartic with $n$ nodes. In these cases the dual graphs correspond
respectively to the diagrams:

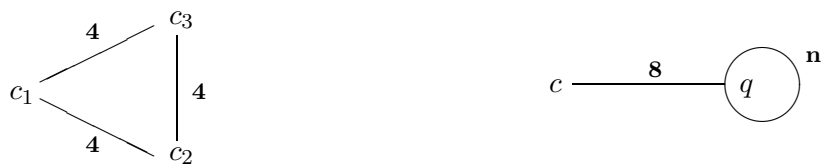
\vspace{0,5cm}

\[
\begin{picture}(0,0)\setlength{\unitlength}{1cm}
\put(-5,-0.2){\line(2,-1){1.5}}\put(-3.2,-0.75){\line(0,3){1.4}}\put(-3.5,+0.77){\line(-2,-1){1.5}}
\put(-5.4,-0.2){$c_1$}\put(-3.3,-1){$c_2$}\put(-3.3,0.8){$c_3$}
\put(-4.4,-0.9){\Small\bf 4}\put(-3,-0.2){\Small\bf
4}\put(-4.4,0.6){\Small\bf 4}
\put(2,0){\line(2,0){2}}\put(4.5,0){\circle{1}}
\put(1.7,-0.1){$c$}\put(4.2,-0.1){$q$}\put(3,0.1){\Small\bf
8}\put(5.07,0.3){\Small\bf n}
\end{picture}\vspace{1cm}
\]
where $c_1,c_2,c_3,c$ are conics, $q$ is a quartic and the bold
number above the edges stands for the number of edges connecting
two vertices. The circle in the second graph means the possible
existence in $q$ of $n$ nodes.

\vspace{0,1cm}

In both cases, to each vertex of the dual graph of $C$ converges
an even number of edges. Moreover, the dual graph of $C$ has first
Betti number greater than zero. (a.1) (with $\tilde C=C$) implies
that there exists an odd theta-characteristic $\theta$ on $C$ such
that $\theta^2=\omega_C$. If $C$ is reducible but it is not as in
Cases 1 and 2, then the dual graph of $C$ is not even. Indeed, $C$
would contain at least a line and such a line would intersect the
union of the other irreducible components of $C$ in an odd number
of points. By the results in \cite{C} and \cite{CC}, there are no
sheaves $L$ on $C$ such that $(C,L)$ is a stable spin curve with
stable model $C$.\end{proof}


The following corollary is a trivial consequence of the techniques
described in the previous proof.


\begin{cor}\label{cor:4.2} Let $X$ be a cubic 4-fold containing
a plane $P$ associated to the pair $(C,\theta)$, where $C$ is a
reduced nodal plane sextic with general irreducible components and
$\theta$ is an odd theta-characteristic with $h^0(C,\theta)=1$.
Then $X$ is smooth if and only if not all the irreducible
components of the total normalization $\tilde C$ of $C$ are
rational.\end{cor}




\begin{remark}\label{rmk:4.1} Assume that a cubic 4-fold $X$ containing a plane $P$ is associated to a reduced nodal plane sextic $C$.
Let $X$ contain ten couples of planes
\[
(P_{1,1},P_{1,2}),\ldots,(P_{10,1},P_{10,2})
\]
which corresponds to ten different fibers of $\pi_{P,C}$. By
Proposition \ref{prop:3.6} and Theorem \ref{thm:3.1}, $P_{i,j}\cap
P_{k,h}$ consists of one point (for $i\neq k$). Then $C$ has one
of the following configurations:

\begin{tabular}{ll}
(a) 6 lines;\hspace{4cm} & (e) 3 lines and 1 cubic with 1 node;\\
(b) 3 conics; & (f) 1 line and a quintic with 5 or 6 nodes;\\
(c) 2 conics and 2 lines; & (g) 1 smooth cubic and 1 cubic with 1
node;\\
(d) 1 conic and 4 lines; & (h) 1 quartic with 1 or 2 nodes and 2
lines.
\end{tabular}

\noindent Moreover, if $X$ is smooth then its discriminant curve
is of type (f) (and it has 5 nodes), (g) or (h) (where the quartic
has just 1 node). Indeed, due to Proposition \ref{prop:3.6},
Theorem \ref{thm:3.1} and Proposition \ref{thm:4.1}, we just need
to prove that (a)--(h) are the only curves which admit a
theta-characteristic $\theta$ such that $h^0(\theta)=1$ and
$\#S_\theta=10$. Since $X$ contains ten couples of planes, $C$
cannot be

(a.1) an irreducible sextic with $m$
nodes ($m<10$);

(b.1) two smooth cubics;

(c.1) a line and a quintic with $m$ nodes ($m<5$);

(d.1) a smooth quadric and a conic,

\noindent because $C$ must have at least ten nodes. We can exclude
the case of an irreducible sextic with ten nodes since the total
normalization of such a curve is a smooth rational curve.

Looking at the dual graph of the remaining possible configurations
for nodal plane sextics, the only ones which are even (i.e. each
vertex has an even number of edges), once we take away ten of their
edges, are those that correspond to sextics as in (a)--(h).

In particular, we get examples of cubic 4-folds containing 10
planes meeting each other in one point and whose discriminant
curve is nodal (see also Example \ref{ex:4.2}(ii)). Theorem
\ref{thm:3.1} implies that some of them can be smooth. Observe
that cubic 4-folds of this type are closely related to the problem
of describing the moduli space of Enriques surfaces (for a more
precise discussion see \cite{Do}).\end{remark}


Now we would like to analyze some explicit examples of singular
cubic 4-folds containing a plane. These will show that the bounds
for the number of points in $\Sing(X)$ given by Theorem
\ref{thm:3.1} are optimal.


\begin{ex}\label{ex:4.2} Let $C$ be a reducible nodal plane sextic
in $\mathbb{P}^5$ and let $x_1,x_2,x_3,u_1,u_2,u_3$ be homogeneous
coordinates in $\mathbb{P}^5$ such that the equations of $\Pi(C)$
are $u_1=u_2=u_3=0$.

(i) Assume moreover that the equation of $C$ in $\Pi(C)$ is the
determinant of the following matrix:
\[
M_1:=\left(
\begin{array}{cccc}
0 & x_1 & x_2 & 0 \\
x_1 & 0 & x_3 & 0 \\
x_2 & x_3 & 0 & 0 \\
0 & 0 & 0 & f
\end{array}
\right),
\]
where the polynomial $f$ defines a degree 3 plane curve $C_1$
meeting the cubic $C_2$, whose equation is $x_1x_2x_3=0$, in nine
distinct points $p_1,\ldots, p_9$ not coinciding with the nodes of
$C_2$. In this case $\tilde S_\theta=\{p_1,\ldots, p_9\}$.

Let $X$ the cubic 4-fold whose equation is
$2x_1u_1u_2+2x_2u_1u_3+2x_3u_2u_3+f=0$. $X$ contains the plane $P$
with equations $x_1=x_2=x_3=0$ and $C$ is a discriminant curve of
$(X,P)$ (by Definition \ref{def:assdisc}). It is very easy to
verify that
\[
B(X,P,C)=\{(0:0:0:1:0:0),(0:0:0:0:1:0),(0:0:0:0:0:1)\}\subseteq P.
\]
Moreover, Proposition \ref{thm:4.1} and Theorem \ref{thm:3.1}
imply that
$\#\Sing(X)=\#B(X,P,C)+\#\Sing(C_1)=3+\#\Sing(C_1)=\#S_C+3$.

(ii) Consider the case when $C$ is a plane curve which is the
union of six lines $l_1,\ldots,l_6$ in general position in
$\Pi(C)\subset\mathbb{P}^5$. Assume that the equation of $C$ in
$\Pi(C)$ is the determinant of the matrix
\[
M_2:=\left(
\begin{array}{cccc}
l_1 & 0 & 0 & 0 \\
0 & l_2 & 0 & 0 \\
0 & 0 & l_3 & 0 \\
0 & 0 & 0 & l_4l_5l_6
\end{array}
\right).
\]
Once more, consider the cubic 4-fold $X$ whose equation is
$l_1u_1^2+l_2u_2^2+l_3u_3^2+l_4l_5l_6=0$ and containing the plane
$P$ with equations $x_1=x_2=x_3=0$. As in the previous case,
$(X,P)$ is associated to $C$ while the base locus $B(X,P,C)$ is
empty. Let $\theta$ be the odd theta-characteristic on $C$ given
by the matrix $M_2$. The set $\tilde S_\theta$ contains 12 of the
15 nodes of $C$. By Theorem \ref{thm:3.1},
\[
\Sing(X)=\Sing(Q_{p_1})\cup\Sing(Q_{p_2})\cup\Sing(Q_{p_3}),
\]
where $\Sing(C)=\tilde S_\theta\cup\{p_1,p_2,p_3\}$ and
$p_1=l_4\cap l_5$, $p_2=l_4\cap l_6$, $p_3=l_5\cap l_6$. Hence
$\#\Sing(X)=\#S_C$.\end{ex}


\begin{ex}\label{ex:4.3} Let $C$ be a plane sextic in $\mathbb{P}^5$ which is the union of a general
quartic $C_1$ and of two lines in general position or of a general
quintic $C_2$ and a line in general position. We choose homogeneous
coordinates $x_1,x_2,x_3,u_1,u_2,u_3$ in $\mathbb{P}^5$ such that
the equations of $\Pi(C)$ are $u_1=u_2=u_3=0$. Since the curves
$C_1$ and $C_2$ are supposed to be general, there are two
theta-characteristics $\theta_1$ and $\theta_2$ such that
\[
h^0(C_1,\theta_1)=1 \qquad \text{and} \qquad h^0(C_2,\theta_2)=1.
\]
By Corollary 4.2 in \cite{B1}, we get the following two matrices
whose determinants are the equations of $C_1$ and $C_2$:
\[
M_1:=\left(\begin{array}{cc} l_{11} & q_1 \\ q_1 & f
\end{array}\right)  \qquad  M_2:=\left(\begin{array}{ccc} l_{11} & l_{12} & q_1 \\ l_{12} & l_{22} & q_2 \\ q_1 & q_2 & f
\end{array}\right).
\]
In these two cases, the sextic $C$ has equation in $\Pi(C)$
described by the determinant of the following two matrices:
\[
N_1:=\left(\begin{array}{ccc} l_1 & 0  & 0\\ 0 & l_2 & 0\\ 0 & 0 &
M_1 \end{array}\right)
\qquad
N_2:=\left(\begin{array}{cc} l_1 & 0 \\ 0 &
M_2\end{array}\right),
\]
where $l_1$ and $l_2$ are the equations of the lines. The
corresponding cubic 4-folds are singular if $C_1$ and $C_2$ have
at least one node.

More explicitly, the matrix corresponding to a sextic which is
union of the Fermat's quartic $x_1^4+x_2^4+x_3^4=0$ and of the two
lines $x_1=0$ and $x_1+x_2=0$ is given by:
\[
M:=\left(
\begin{array}{cccc}
-x_1 & 0 & 0 & 0 \\
0 & x_1+x_2 & 0 & 0 \\
0 & 0 & -(x_1-\omega x_2) & x_3^2 \\
0 & 0 & x_3^2 & (x_1+\omega x_2)(x_1^2+ix_2^2)
\end{array}
\right),
\]
where $i,\omega\in\mathbb{C}$ are such that $i^2=-1$ and
$\omega^2=-i$.\end{ex}


The following theorem proves, in particular, that there is a
smooth rational cubic 4-fold $X$ containing a plane $P$ associated
to a reducible nodal plane sextic. We write $\mathrm{Mat}(3)$ for
the algebra of $3\times 3$ matrices with complex coefficients.


\begin{prop}\label{thm:4.4} Given 3 lines in $\mathbb{P}^5$ meeting in three distinct points, there exists a
family of smooth rational cubic 4-folds $X$ containing a plane $P$
such that a discriminant curve of the pair $(X,P)$ is a reduced
nodal plane sextic containing these three lines. This family is
parametrized by points in a non-empty open subset of
$\mathrm{Mat}(3)$. Moreover, the group $\mathrm{NS}_2(X)$ contains
25 distinct classes corresponding to planes in $X$ and
$\mathrm{rkNS}_2(X)\geq 14$.

In particular, there exists a smooth rational cubic 4-fold
containing a plane associated to reduced nodal plane
sextic.\end{prop}


\begin{proof} Let $l_1,l_2,l_3\subset\mathbb{P}^5$ be three lines meeting each
other in three distinct points. Fix homogeneous coordinates
$x_1,x_2,x_3,u_1,u_2,u_3$ in $\mathbb{P}^5$ such that the equations
of $l_1$, $l_2$ and $l_3$ are
\[
l_1: u_1=u_2=u_3=x_1=0\;\;\;\;\;\;\; l_2: u_1=u_2=u_3=x_2
=0\;\;\;\;\;\;\; l_3: u_1=u_2=u_3=x_3=0.
\]
The lines $l_1$, $l_2$ and $l_3$ are contained in the plane $\Pi$
whose equations are $u_1=u_2=u_3=0$.

If $P$ is the plane described by the equations $x_1=x_2=x_3=0$,
given a matrix $A:=(a_{ij})\in \mathrm{Mat}(3)$, the equations of
a plane $P_A$ such that $P\cap P_A=\emptyset$ can be obtained as
the zero-locus of
\begin{align*}
u_1-a_{11}x_1-a_{12}x_2-a_{13}x_3=0 \\
u_1-a_{21}x_1-a_{22}x_2-a_{23}x_3=0 \\
u_1-a_{31}x_1-a_{32}x_2-a_{33}x_3=0.
\end{align*}

Let $C'_A$ be the plane cubic in $\Pi$ defined by the polynomial
\[
f_A:=\sum_{i=1}^3(a_{i1}x_1+a_{i2}x_2+a_{i3}x_3)^2x_i.
\]
If we require that $C'_A$ is smooth and that $C'_A\cap(l_1\cup
l_2\cup l_3)\neq\{(1:0:0),(0:1:0),(0:0:1)\}$, we impose that $A$
belongs to an open subset $U\subseteq \mathrm{Mat}(3)$. Given
$A\in U$, we get a nodal reduced plane sextic $C_A:=C'_A\cup
l_1\cup l_2\cup l_3$. By Proposition \ref{prop:2.4}, the matrix
\[
M_A:=\left(
\begin{array}{cccc}
x_1 & 0 & 0 & 0 \\
0 & x_2 & 0 & 0 \\
0 & 0 & x_3 & 0 \\
0 & 0 & 0 & -f_A
\end{array}
\right)
\]
determines a theta-characteristic $\theta_A$ on $C_A$ such that
$h^0(C_A,\theta_A)=1$. Consider the cubic 4-fold $X_A$ whose
equation is
\[
F_A:=x_1u_1^2+x_2u_2^2+x_3u_3^2-f_A.
\]
Obviously $P,P_A\subset X_A$ and $(X_A,P)$ is associated to
$(C_A,\theta_A)$. As it was shown in \cite{T1} and \cite{T2}, the
planes $P_A$ gives a section for the projection
$\pi_{P,C_A}:X_A\dashrightarrow\Pi(C_A)$ and $X_A$ is rational. It
is very easy to see that $\tilde S_{\theta_A}=\Sing(C_A)$ and that
$B(X_A,P,C_A)=\emptyset$. Hence Theorem \ref{thm:3.1} implies that
$X_A$ is smooth.

To show that the open subset $U$ is non-empty, let $C'$ be the
Fermat plane cubic in $\Pi$ whose equation in $\Pi$ is
$f=x_1^3+x_2^3+x_3^3$. Clearly, the three lines $l_1$, $l_2$ and
$l_3$ meet the cubic in nine points distinct from the three
intersection points $(1:0:0)$, $(0:1:0)$ and $(0:0:1)$. By simple
calculations we see that the plane $P'$ described by the equations
\[
u_1-x_1=u_2-x_2=u_3-x_3=0
\]
is contained in the cubic 4-fold $X$ whose equation in
$\mathbb{P}^5$ is $x_1u_1^2+x_2u_2^2+x_3u_3^2-f=0$. Moreover,
$P\cap P'=\emptyset$. Hence $X$ is rational. This implies that the
matrix $\mathrm{Id}\in \mathrm{Mat}(3)$ is in $U$.

Given $A\in U$, we have the equalities $S_{\theta_A}=\Sing(C_A)$
and $\#S_{\theta_A}=12$. Theorem \ref{thm:3.1} gives the desired
estimate about the number of distinct planes in
$\mathrm{NS}_2(X_A)$ and the rank of
$\mathrm{NS}_2(X_A)$.\end{proof}

\medskip

{\small\noindent{\bf Acknowledgements.} During the preparation of
this paper the author was partially supported by the MIUR of the
Italian Government in the framework of the National Research
Project ``Geometry on Algebraic Varieties'' (Cofin 2002). This
work was started during the summer school PRAGMATIC 2003 which
took place in Catania. The author would like to thank Sandro Verra
and Igor Dolgachev for their encouragement and the organizers for
their hospitality. It is a pleasure to thank Bert van Geemen for
helpful discussions and Sandro Verra for his suggestions and for
pointing out some simplifications in the proof of Theorem
\ref{thm:3.1}. Thanks also to Rogier Swierstra for carefully
reading these pages and to Marco Pacini for his suggestions about
some references related to theta-characteristics on curves.}


\end{document}